
\magnification=1200
\hfuzz=3pt
\overfullrule=0mm


\font\tensymb=msam9
\font\fivesymb=msam5 at 5pt
\font\sevensymb=msam7  at 7pt
\newfam\symbfam
\scriptscriptfont\symbfam=\fivesymb
\textfont\symbfam=\tensymb
\scriptfont\symbfam=\sevensymb

\font\titlefont=cmbx10 at 15pt


\font\refttfont=cmtt10 at 9pt

\font\sc=cmcsc10 \rm


\def\Hom{{\rm Hom}}
\def\Inf{{\rm Inf}}
\def\Ker{{\rm Ker}}

\def\Fq{{\bf F}}
\def\ZZ{{\bf Z}}

\def\N{{\rm N}}
\def\and{\quad\hbox{and}\quad}

\def\Pr{\noindent {\sc Proof.--- }}

\def\cqfd{ {\sevensymb {\char 3}}}

\def\hfl#1#2{\smash{\mathop{\hbox to 6mm{\rightarrowfill}}
\limits^{\scriptstyle#1}_{\scriptstyle#2}}}

\def\vfl#1#2{\llap{$\scriptstyle #1$}\left\downarrow
\vbox to 3mm{}\right.\rlap{$\scriptstyle #2$}}


\null
\noindent

\vskip 20pt

\noindent
\centerline{\titlefont
Explicit elements of norm one for cyclic groups}

\vskip 30pt

\centerline{\sc Eli Aljadeff and Christian Kassel}

\footnote{}{Supported by TMR-Grant ERB FMRX-CT97-0100 of the European Commission}

\vskip 15pt
\bigskip

\noindent
Let $R$ be a ring (with unit element~$1$) and $G$ a finite group
acting on~$R$ by ring automorphisms. For any subgroup~$U$ of~$G$ define
the norm map (sometimes called the trace map) $\N_U : R\to R^U$ by
$$\N_U(x)  = \sum_{g\in U}\, g(x).$$
Ginosar and the first author reduced the surjectivity
of the norm map $\N_G$ for a group~$G$ to the surjectivity of the 
norm maps for its elementary abelian subgroups; more precisely,
they proved that $\N_G: R \to R^G$ is surjective if and only if
$\N_U: R\to R^U$ is surjective for every elementary abelian subgroup~$U$
of~$G$ (see [2, Theorem~1]). 

The $R^U$-linearity of $\N_U$ implies that it is surjective if and
only there exists an element $x_U\in R$ such that $\N_U(x_U) = 1$.
Suppose we have such an element~$x_U$ for every elementary abelian
subgroup~$U$ of~$G$. Then by the result mentioned above there is a
``global" element $x_G\in R$ such that~$\N_G(x_G) = 1$. 
Using this last statement, Shelah observed (see [2,~Proposition~6]) 
that there exists a formula in which $x_G$ is a finite sum of the form
$$x_G = \sum\,
a \, \sigma_{i_{1}}(x_{_{U_{j_{1}}}})\, \sigma_{i_{2}}(x_{_{U_{j_{2}}}})\,
\cdots\, \sigma_{i_{r}}(x_{_{U_{j_{r}}}}),$$ 
where $a\in\ZZ$ and $\sigma_{i_{1}}, \sigma_{i_{2}}, \ldots,
\sigma_{i_{r}}\in G$.

Using a tensor induction argument, the first author found such a formula
in case the ring~$R$ is commutative (see [1, Theorem~2.1]). 
When $R$ is not commutative, the only formulas known so far hold in the following
two cases:

(a) $G$ is an abelian $2$-group and $R$ is an algebra over the
field ~$\Fq_2$ ({\it cf}.~[2, Section~2]),

(b) $G = \ZZ/4$ and $R$ is any ring (see Formula~(2) below).

The aim of this article is to show how to find an explicit formula of the above form
in the case of a cyclic $p$-group~$G$ acting on an arbitrary (non-necessarily
commutative) ring~$R$. 
More precisely, in the theorem below, we express a norm one element for a cyclic group
$G$ of order~$p^n$ in terms of a norm one element for a
subgroup of order~$\geq p^{n/2}$. This allows to find by induction
a norm one element for a cyclic $p$-group in terms of a
norm one element for its unique elementary subgroup ($\cong \ZZ/p$).

The theorem appears as the main result of Section~1.
Its proof is given in Section~2.
In Section~3 we give some cohomological explanations for the proof.

\bigskip\goodbreak
\noindent
{\bf 1.\ Statement of results}

\medskip
\noindent
We fix a prime number~$p$, and integers $n$ and $k$ verifying $n\geq 2$
and $1\leq k \leq n/2$.
Consider the cyclic group $G = \ZZ/p^n$ of order~$p^n$ with a
generator~$\sigma$, and the cyclic subgroup~$U$ of order~$p^{n-k}$,
generated by~$\sigma^{p^k}$.
Our main result is the following.

\medskip\goodbreak
\noindent
{\sc Theorem}.---
{\it Let $x\in R$ satisfy $\N_U(x) = 1$. Define $z$ and $a\in R$ by
$$z = p^{n-2k}\,
\bigl(1 + \sigma + \sigma^2 + \cdots + \sigma^{p^k -1}\bigr) (x) -1$$
and
$$a = p^{n-2k} x + (1-\sigma)
\biggl(\;  \sum_{i=1}^{p^{n-k} -1}\,
\Bigl( 1 + \sigma^{p^k} + \sigma^{2p^k} + \cdots + \sigma^{(i-1)p^k}
\Bigr)
\bigl( x\sigma^{-ip^k}(z)
\bigr)
\biggr). \eqno (1)$$
Then $\N_G(ax) = 1$.
}

\medskip\medskip
When $G = \ZZ/p^2$, the theorem can be reformulated as follows.

\medskip\medskip\goodbreak
\noindent
{\sc Corollary}.---
{\it Let $\sigma$ be an automorphism of order~$p^2$ of a ring~$R$,
and $x\in R$ satisfying
$( 1 + \sigma^{p} + \sigma^{2p} + \cdots + \sigma^{(p-1)p})(x) = 1.$
Define
$z =
\bigl(1 + \sigma + \sigma^2 + \cdots + \sigma^{p -1}\bigr) (x) -1$
and
$$a = x + (1-\sigma)
\biggl(\;  \sum_{i=1}^{p -1}\,
\Bigl( 1 + \sigma^{p} + \sigma^{2p} + \cdots + \sigma^{(i-1)p}
\Bigr)
\bigl( x\sigma^{-ip}(z) \bigr) \biggr).$$
Then
$$ (1 + \sigma + \sigma^2 + \cdots + \sigma^{p^2-1})(ax) = 1.$$
}

\medskip
For $p = 2, 3$, the corollary yields the following explicit formulas.
When $p=2$, starting from $x\in R$ such that $x+\sigma^2(x) = 1$, we
obtain the norm one element
$$\eqalign{
ax & =  \sigma(x)x  - \sigma(x) x^2 + x \sigma^2(x)x + x\sigma^3(x)x
- \sigma(x)\sigma^3(x)x\cr
& = 2x^2 - x^3 - x\sigma(x)x  - \sigma(x)x^2 +\sigma(x)^2x.\cr 
}$$
In this case, P\'eter P.~P\'alfy had shown the first author the following
simpler formula:
$$x\sigma(x)x + x\sigma(x) - x^2\sigma(x).  \eqno (2)$$

When $p=3$, starting now from $x\in R$ such that
$x+\sigma^3(x) + \sigma^6(x) = 1$, we obtain the norm one element
$$\eqalign{
ax = &
- x^2 + 2 \sigma(x)x - \sigma^3(x)x + \sigma^4(x)x \cr
&{} + x\sigma^3(x)x + x\sigma^4(x)x + x\sigma^5(x)x + x\sigma^6(x)x
+ x\sigma^7(x)x + x\sigma^8(x)x\cr
 & - \sigma(x)\sigma^4(x)x - \sigma(x)\sigma^5(x)x
- \sigma(x)\sigma^6(x)x - \sigma(x)\sigma^7(x)x -
\sigma(x)\sigma^8(x)x - \sigma(x) x^2 \cr
& + \sigma^3(x)\sigma^6(x)x + \sigma^3(x)\sigma^7(x)x +
\sigma^3(x)\sigma^8(x)x\cr
& - \sigma^4(x)\sigma^7(x)x -  \sigma^4(x)\sigma^8(x)x
- \sigma^4(x)x^2.\cr
}$$

\medskip
Using the theorem repeatedly, we can find an explicit expression
for an element~$x_G\in R$ with $\N_G(x_G)=1$ as a noncommutative
polynomial in the variables $\sigma^i(x_E)$ ($0\leq i < p^n$),
where $x_E\in R$ satisfies $\N_E(x_E) = 1$ for the unique
subgroup~$E\cong \ZZ/p$ of~$G$.
Formula~(1) also gives an upper bound for the number of monomials
appearing in this polynomial.
Indeed, the number of monomials for~$a$ in~(1) is
$\leq p^{n-k}(p^{n-k} -1)(p^{k}+1) +1$.
We obtain a crude upper bound by setting $k=1$ and by summing from
$2$ to~$n$. The upper bound we get in this way is
$${p(p+1)(p^{n-1} - 1) (p^n - 1) \over p^2 - 1} + n - 1 
\sim {p+1 \over p^2 - 1}\, p^{2n}$$
when $n$ becomes large.
By taking~$k$ as the largest integer~$\leq n/2$, we
get a smaller upper bound, which is equivalent
to~$p^{1/2}\, p^{3n/2}$~($n \gg 0$).

\bigskip\goodbreak
\noindent
{\bf 2.\ Proof of the theorem}

\medskip
\noindent
We shall need a special case of the following lemma,
directly inspired from Proposition~1.3 of [3, Chap.~XII].

\medskip\goodbreak
\noindent
{\sc Lemma~1}.---
{\it Let $U$ be a finite group acting by ring automorphisms on a
ring~$R$. If there exists an element $x\in R$ such that $\N_U (x) = 1$,
then every element $z\in R$ such that $\N_U (z)= 0$ can be written as
$$z = \sum_{g\in U}\, (g-1) \bigl( xg^{-1}(z) \bigr) .$$
}

\medskip
\Pr
We have
$$\eqalign{
\sum_{g\in U}\, (g-1) \bigl( xg^{-1}(z) \bigr)
& = \sum_{g\in U}\,  g(x) \, g\bigl( g^{-1}(z) \bigr)
- \sum_{g\in U}\, xg^{-1}(z)  \cr
& = \N_U(x)\, z  - x\, \N_U(z) = z.\cr
}$$
\hfill\cqfd
\medskip

Let us apply Lemma~1 to the case when $U$ is a
cyclic group of order~$r$. If we denote a generator of~$U$ by~$t$, then
every element $z\in R$ such that
$\N_U (z) = 0$ is of the form $z = (t-1) (w)$, where
$$w= \sum_{i=1}^{r-1} \, (1  + t + t^2 + \cdots + t^{i-1})
\bigl( xt^{-i}(z) \bigr) . \eqno (3)$$

We now start the proof of the theorem.
Recall that $G$ is a cyclic group of order~$p^n$, generated
by~$\sigma$, and $U$ is the subgroup generated
by~$\sigma^{p^k}$, where $1\leq k \leq n/2$.

Consider the group $B= \Hom(\ZZ[G],R)$ of $\ZZ$-linear maps from the
group ring~$\ZZ[G]$ to~$R$, and equip it with the $G$-action given by
$(g\varphi)(s)  = \varphi(sg)$
for all $g,s\in G$ and $\varphi \in B$.
We embed $R$ into $B$ by considering an element $x\in R$ as
the element $\varphi_x \in B$ determined by $\varphi_x(g) = g(x)$
for all $g\in G$. The embedding $R\subset B$ preserves the $G$-action.

Define $\varphi \in B$ by 
$$\varphi(g) = \cases{ 1 & if $g =\sigma^{ip^k}$ $( 0 \leq i < p^{n-k})$,
\cr
\noalign{\smallskip}
$0$ & otherwise.\cr} 
$$ 
It is clear that $\varphi$ is invariant
under the action of the subgroup~$U$. Since the action of $G$ on
the group~$B^U$ of $U$-invariant elements of~$B$ factors through
an action of the quotient cyclic group~$G/U$, we may consider the
action of 
$$\N_{G/U} = 1 +  \sigma + \sigma^2 + \cdots +
\sigma^{p^k-1}$$ 
on~$B^U$. We clearly have 
$$\N_{G/U}(\varphi) =
\varphi_1 \and \N_G(\varphi) = p^{n-k}\, \varphi_1 . \eqno (4)$$
On the other hand, $\N_U(x) = 1$ implies $\N_G(x) = p^k$.
Therefore 
$$\N_G(\varphi - p^{n-2k} \, \varphi_x) = 0. \eqno (5)$$
Define $\psi \in B$ inductively by $\psi(\sigma^0) = 0$ and for
$1\leq i < p^n$ by 
$$\psi(\sigma^{i}) = \psi(\sigma^{i-1}) -
\varphi(\sigma^{i-1}) + p^{n-2k} \, \sigma^{i-1}(x).$$ 
It follows from the definition of~$\psi$ and from (5) that 
$$\varphi = (1-\sigma) (\psi) + p^{n-2k} \varphi_x. \eqno (6)$$

\medskip\goodbreak
\noindent
{\sc Lemma~2}.---
{\it The element $\psi \in B$ is $U$-invariant modulo~$R$.
}

\medskip
\Pr
By (4) and (6) we have
$$\eqalign{
\varphi_1  & = \N_{G/U}(\varphi)  \equiv
\N_{G/U}\bigl( (1-\sigma) (\psi)\bigr) \cr
& \equiv \bigl((1 +  \sigma + \sigma^2 + \cdots + \sigma^{p^k-1})
(1-\sigma)\bigr) (\psi) \cr
& \equiv (1 -  \sigma^{p^k})(\psi) \quad\hbox{\rm modulo} \, R.
}$$
We conclude by observing that $\varphi_1 \equiv \varphi_0 \equiv 0$
modulo~$R$.
\hfill\cqfd
\medskip

It follows from Lemma~2 that there exists $z\in R$ such that
$$\bigl( \sigma^{p^k}-1\bigr)(\psi) = \varphi_z. \eqno (7)$$
In order to find~$z$, it suffices to apply both sides of~(7)
to the unit element of~$G$; we thus obtain
$z = \psi(\sigma^{p^k}) - \psi(\sigma^0)$, which, by definition of~$\psi$,
equals
$$z = p^{n-2k}\, \bigl(1 +  \sigma + \sigma^2 + \cdots +
\sigma^{p^k-1}\bigr)(x) - 1.\eqno (8)$$

We claim that the element $z\in R$ is killed by~$\N_U$. Indeed, by~(7),
$$\N_U( z) = \N_U\bigl((\sigma^{p^k} - 1)(\psi)\bigr)
= \bigl(\sigma^{p^n} - 1\bigr)(\psi) = 0.$$
Applying Formula~(3), we have $z = (\sigma^{p^k}-1) (w)$, where
$$w= \sum_{i=1}^{p^{n-k}-1} \,
\bigl(1  + \sigma^{p^k} + \sigma^{2p^k} + \cdots + \sigma^{(i-1)p^k}
\bigr) \bigl( x\sigma^{-ip^k}(z) \bigr) . \eqno (9)$$
Formulas (7) and (9) imply
$$\bigl( \sigma^{p^k} -1 \bigr) (\psi - \varphi_w) = 0. \eqno (10)$$
In other words, $\psi - \varphi_w $ is $U$-invariant.

\medskip\medskip\goodbreak
\noindent
{\sc Lemma~3}.---
{\it The element $a = p^{n-2k}\, x + (1 - \sigma)(w) \in R$ is
$U$-invariant and $\N_{G/U}(a) = 1$.
}

\medskip
\Pr
By definition of $a$ and $w$, and by~(6), we have
$$\varphi_a = \varphi - (1-\sigma)(\psi - \varphi_w),$$
which is $U$-invariant in view of~(10). By (4) and (10), we obtain
$$\eqalign{
\N_{G/U}(a)
& = \N_{G/U}(\varphi)
- \N_{G/U}\bigl( (1-\sigma)(\psi -\varphi_w)\bigr) \cr
& = \varphi_1 - (1-\sigma^{p^k})(\psi - \varphi_w)  = 1. \cr }$$
\hfill\cqfd

\medskip
We can now complete the proof of the theorem. First observe that the
elements $a\in R^U$ and $z\in R$ of the theorem are exactly the ones introduced
in this section. So it is enough to check that $\N_G(ax) = 1$.
Indeed, using the $R^U$-linearity of~$\N_U$, Lemma~3, and
$\N_U(x) = 1$, we have
$$\eqalign{
\N_G(ax) & = \sum_{g\in G}\, g(ax)
= \sum_{t\in G/U}\, t\Bigl( \sum_{s\in U}\,  s(ax) \Bigr) \cr
& = \sum_{t\in G/U} \, t\bigl( \N_U(ax) \bigr)
= \sum_{t\in G/U} \, t\bigl( a\N_U(x)\, \bigr)\cr
& = \sum_{t\in G/U} \, t(a) = \N_{G/U}(a) = 1. \cr
}$$
\hfill\cqfd

\bigskip\goodbreak
\noindent
{\bf 3.\ Some cohomological considerations}

\medskip
\noindent
In the proof of the theorem as given in Section~2, the computations take place
in the co-induced module $B = \Hom(\ZZ[G],R)$, and the elements
$\varphi$, $\psi\in B$ play a central r\^ole. This can be explained
through the following cohomological considerations.

We first claim that the existence of $x\in R$ such that $\N_U(x) = 1$ implies
the vanishing of the cohomology of the group~$U$ with coefficients in~$R$
in positive degrees:
$H^i(U,R) = 0$ for all~$i>0$.
Indeed, since $U$ is cyclic with generator~$\sigma^{p^k}$, we have
$$H^1(U,R) = \Ker(\N_U: R\to R)/(\sigma^{p^k} - 1)(R),$$
which is zero by Lemma~1.
The surjectivity of $\N_U : R\to R^U$ implies the vanishing of
$$H^2(U,R) = R^U/\N_U(R).$$
Then the claim follows from the periodicity of the cohomology of
cyclic groups.

In view of~[4, Section~VII.6] (or of Hochschild-Serre's spectral sequence),
the vanishing of $H^i(U,R) = 0$ for~$i>0$ implies that the inflation maps
$$\Inf : H^*(G/U,R^U) \to H^*(G,R)$$
are isomorphisms.

Now consider the short exact sequence of $\ZZ[G]$-modules
$$0 \to R \to B \to C \to 0, \eqno (11)$$
where $C = B/R$.
Applying $H^*(U,-)$ to~(11), we obtain a sequence
of $\ZZ[G/U]$-modules
$$0 \to R^U \to B^U \to C^U \to 0, \eqno (12)$$
which is exact because of the vanishing of~$H^1(U,R)$.
Observe that $B^U \cong \Hom(\ZZ[G/U],R)$ is a co-induced module
for~$G/U$.

The naturality of the inflation maps gives rise to the commutative square
$$\matrix{
H^1(G/U,C^U) & \hfl{\Inf}{} & H^1(G,C)\cr
\noalign{\smallskip}
\vfl{\delta}{} & & \vfl{\delta}{} \cr
\noalign{\smallskip}
H^2(G/U,R^U) & \hfl{\Inf}{} & H^2(G,R)\cr
} \eqno (13)$$
where the vertical maps $\delta$ are connecting maps
for the short exact sequences~(11) and~(12).
We claim that all maps in the square~(13) are isomorphisms. We have already
proved this for the lower inflation map. The connecting maps~$\delta$
are isomorphisms because co-induced modules are cohomologically trivial.
It follows that the upper inflation map is an isomorphism as well.
Moreover, by [2, Theorem~1], the surjectivity of $\N_U: R\to R^U$
implies the surjectivity of $\N_G: R\to R^G$. Therefore,
$H^2(G,R) = R^G/\N_G(R) = 0$,
which implies the vanishing of all cohomology groups in~(13).

The central r\^ole played by the element~$\varphi \in B^U$
in the proof of the theorem follows from the following two facts:

(i) If $\bar{\varphi}$ denotes the class of $\varphi$ in~$C$, then
by~(4) it induces an element
$$[\bar{\varphi}] \in H^1(G/U,C^U)
= \Ker(\N_{G/U} : C^U\to C^U)/(\sigma-1)(C^U).$$

(ii) The image $\delta([\bar{\varphi}])$ of $[\bar{\varphi}]$ under
the connecting map for the short exact sequence~(12) is computed as follows:
lift $\bar{\varphi}$ to $\varphi\in B^U$ and apply~$\N_{G/U}$.
By~(4) again, we obtain
$$\delta([\bar{\varphi}]) = 1 \in H^2(G/U,R^U) = R^G/\N_{G/U}(R^U).$$

The existence of $\psi$ satisfying~(6) follows from the vanishing of
$$H^1(G,C) = \Ker(\N_G : C\to C)/(\sigma-1)(C)$$
and Lemma~2 follows from the injectivity of $\Inf : H^1(G/U,C^U) \to
H^1(G,C)$.

\bigskip\bigskip\goodbreak
\centerline{\bf References}
\bigskip

\noindent
[1] {\sc E.~Aljadeff},
{\it On the surjectivity of some trace maps},
Israel J.~Math.\ 86 (1994), 221--232.
\smallskip

\noindent
[2] {\sc E.~Aljadeff, Y.~Ginosar},
{\it Induction from elementary abelian subgroups},
J.~of Algebra 179 (1996), 599--606.
\smallskip

\noindent
[3] {\sc H.~Cartan, S.~Eilenberg},
{\it Homological algebra},
Princeton University Press, Princeton, 1956.
\smallskip

\noindent
[4] {\sc J.-P.\ Serre},
{\it Corps locaux}, Publications de l'Universit\'e de Nancago, Hermann, Paris, 1962
(English translation: {\it Local fields},
Grad.\ Texts in Math.\ 67, Springer-Verlag, New York, Berlin, 1979).
\smallskip

\vskip 30pt
\line{Eli Aljadeff\hfill}
\line{Department of Mathematics \hfill}
\line{Technion - Israel Institute of Technology \hfill}
\line{32000 Haifa, Israel \hfill}
\line{E-mail:  {\refttfont aljadeff@techunix.technion.ac.il}\hfill}
\line{Fax: +972-4-832-4654 \hfill}
\medskip\medskip

\line{Christian Kassel\hfill}
\line{Institut de Recherche Math\'ematique Avanc\'ee \hfill}
\line{C.N.R.S. - Universit\'e  Louis Pasteur\hfill}
\line{7 rue Ren\'e Descartes \hfill}
\line{67084  Strasbourg Cedex, France \hfill}
\line{E-mail:  {\refttfont kassel@math.u-strasbg.fr}\hfill}
\line{Fax: +33 (0)3 88 61 90 69 \hfill}
\line{http:/\hskip -1.5pt/www-irma.u-strasbg.fr/\raise-4pt\hbox{\~{}}kassel/\hfill}

\vskip 25pt
\noindent
(5 December 2000)

\bye